\documentclass[12pt,oneside]{aomamlt2e}
\theoremstyle{plain}
\newtheorem{theorem}{Theorem}[section]
\newtheorem{cor}[theorem]{Corollary}
\newtheorem{prop}[theorem]{Proposition}
\newtheorem{lemma}[theorem]{Lemma}

\theoremstyle{plain}
\newtheorem{de}[theorem]{Definition}
\newtheorem{rem}[theorem]{Remark}
\newtheorem{nota}[theorem]{Notation}
\newtheorem{property}[theorem]{Property}

\numberwithin{equation}{section}

\newcommand{\E}{\mathcal E}
\newcommand{\N}{\mathcal N}
\renewcommand{\L}{\mathcal L}
\newcommand{\V}{\mathcal V}

\newcommand{\Sm}{\it Sm/k}

\renewcommand{\P}{\mathbb P}
\newcommand{\A}{\mathbb A}
\newcommand{\triv}{\mathbf{1}}

\newcommand{\GL}{\mathcal{GL}}

\newcommand \lb{\lbrack}
\newcommand \rb{\rbrack}

\newcommand \id{\text{id}}

\newcommand \spe{\operatorname{Spec}}

\renewcommand \id{\operatorname{id}}

\newcommand \pt{\operatorname{pt}}

\renewcommand \phi\varphi
\renewcommand \Im {\operatorname{Im}}

\def \rra#1{\overset{#1}{\rightarrow}}

\def \ab{\mathcal Ab}
\def \ab2{{\mathbb Z}/{\rm 2\text-}\mathcal Ab}

\def \dia {\diamond}

\newcommand{\mono}{\hookrightarrow}

\renewcommand{\cap}{\smallfrown}
\renewcommand{\cup}{\smallsmile}

\begin{document}

\title{$T$-spectra and Poincar\'e Duality}
\twoauthors{Ivan Panin\footnote{Supported in part by the Ellentuck
Fund, by the RTN-Network HPRN-CT-2002-00287 and by the Presidium
of RAS Program ``Fundamental Research''.}}{Serge Yagunov\footnote
{Both authors were partially supported by RTN Network
HPRN-CT-2002-00287, INTAS, and the Russian Academy of Sciences
research grants from the "Support Fund of National Science" 2001-3
(the first author) and for 2004-5 (the second one).}}

\begin{abstract}

Frank Adams introduced the notion of a complex oriented cohomology
theory represented by a commutative ring-spectrum and proved the
Poincar\'e Duality theorem for this general case. In the current
paper we consider oriented cohomology theories on algebraic
varieties represented by  multiplicative symmetric $T$-spectra and
prove the Duality theorem, which
 mimics the result of Adams.
This result is held, in particular, for Motivic Cohomology and
Algebraic Cobordism of Voevodsky.

\end{abstract}

\setcounter {section}{-1}
\section{Introduction}

In certain cases a commutative
ring-spectrum $E$
can be equipped
with
a distinguished element
$c \in E^2(\P^{\infty})$
called a complex orientation of $E$ (see \cite{Ad}).
The pair $(E,c)$ is called a complex oriented ring spectrum.
Given a complex orientation $c$ of $E$,
every smooth complex projective variety $X$
can be equipped with a homological class
$\lb X\rb \in E_{2d}(X)$
called the fundamental class of $X$
(here $d$ stays for the complex
dimension of $X$). This class has the property that the cap-product
$$
\cap\!\!\lb X\rb: E^*(X) \to E_{2d - *}(X)
$$
provides an isomorphism of cohomology
and homology groups of $X$. The isomorphism
is often called the Poincar\'e Duality isomorphism.

From the modern point of view it looks pretty interesting to
obtain an analogue of this result in the context of Algebraic
Geometry. It is reasonable in this case to choose and fix a field
$k$ and consider a symmetric commutative ring $T$-spectrum
$\mathcal A$ in the sense of Voevodsky~\cite{Vo} (for the
concept of symmetric $T$-spectrum see
Jardine~\cite{Ja}). The $T$-spectrum $\mathcal A$
determines bi-graded cohomology and homology theories ($A^{*,*}$
and $A_{*,*}$) on the category of algebraic varieties
~\cite[p.595]{Vo}. (We also assume the spectrum $A$ to be a ring-spectrum {\it i.e.}
be endowed with a multiplication $\mu\colon\mathcal A\wedge\mathcal A\to\mathcal A$,
which induces product structures in (co)homology.)
 In some cases $\mathcal A$ can be equipped with
a distinguished element $\gamma \in A^{2,1}(\P^{\infty})$, which
Morel calls an orientation of $\mathcal A$. Following him the pair
$(\mathcal A,\gamma)$ is called an oriented symmetric commutative
ring $T$-spectrum. The orientation $\gamma$ equips both cohomology
$A^{*,*}$ and homology $A_{*,*}$ with trace
structures~\cite{Pa2,Pi}. The latter means that for every
projective morphism $f: Y \to X$ of $k$-smooth irreducible
varieties with $d=\text{dim}(X)-\text{dim}(Y)$ there are two
operators $f_!: A^{*,*}(Y) \to A^{*+2d,*+d}(X)$ and $f^!:
A_{*,*}(X) \to A_{*-2d,*-d}(Y)$ satisfying a list of natural
properties. Define now a fundamental class of a $k$-smooth
projective equi-dimensional variety $X$ of dimension $d$ as $\lb
X\rb \overset{\text{def}}=\pi^!(1)\in A_{2d,d}(X)$, where
$\pi\colon X\to\pt$ is the structure morphism. Our main result
claims that the map
\begin{equation*}
\cap\lb X\rb\colon A^{*,*}(X)\overset\simeq\to A_{2d-*,d-*}(X)
\end{equation*}
is a grade-preserving isomorphism (Poincar\'e Duality isomorphism).

There are at least two interesting examples of oriented symmetric
commutative ring $T$-spectra. The first one is a symmetric model
$\mathbb {MGL}$ of the algebraic cobordism $T$-spectrum
$\textrm{MGL}$ of Voevodsky~\cite[p. 601]{Vo}. This symmetric
commutative ring $T$-spectrum $\mathbb {MGL}$ together with an
orientation $\gamma \in \mathbb {MGL}^{2,1}(\P^{\infty})$ is
described in details in~\cite[Sect.6.5]{PY}. So that, every smooth
irreducible projective variety $X/k$ of dimension $d$ has the
fundamental class $[X] \in \mathbb {MGL}_{2d,d}(X)$ and the
cap-product with this class
\begin{equation*}
\label{DualityForCobordism}
\cap\lb X\rb\colon \mathbb {MGL}^{*,*}(X)\overset\simeq\to \mathbb {MGL}_{2d-*,d-*}(X)
\end{equation*}
is an isomorphism.

The second example is the Eilenberg--Mac Lane $T$-spectrum
$\textrm{H}$ (it is intrinsically a symmetric $T$-spectrum
representing the motivic cohomology). This $T$-spectrum
$\textrm{H}$ is constructed in~\cite[p.598]{Vo} and we briefly
describe its orientation here. Recall that for a smooth variety
$X/k$ the first Chern class of a line bundle with value in the
motivic cohomology defines a functorial isomorphism
$\textrm{Pic}(X)= \textrm H^{2,1}(X)$. Thus, $\mathbb Z= \textrm
H^{2,1}(\P^{\infty})$ and the class of the line bundle $\mathcal
O(-1)$ over $\P^{\infty}$ is a free generator of $\textrm
H^{2,1}(\P^{\infty})$. This class provides the required
orientation of $\textrm{H}$. Similarly to the case of algebraic
cobordism, one has the fundamental class $[X] \in \textrm
H^\mathcal M_{2d,d}(X)$ in Motivic homology and the isomorphism:
\begin{equation*}
\label{DualityForMotivicCoh}
\cap\lb X\rb\colon \textrm{H}_\mathcal M^{*,*}(X)\overset\simeq\to \textrm{H}^\mathcal M_{2d-*,d-*}(X).
\end{equation*}
To embellish the latter result, let us mention that unlike the topological context in algebro-geometrical case the canonical pairing
$\textrm{H}_\mathcal M^{*,*}(X)\otimes\textrm{H}^\mathcal M_{*,*}(X)\to \textrm{H}_\mathcal M^{*,*}(\pt)$
is generally degenerated even with rational coefficients~\cite{Vo2}.

The paper is organized as follows. Section 1 is devoted to
product structures in extraordinary cohomology and homology
theories. In section 2 we formulate Poincare Duality theorem and
derive it from two projection formulae, which are proven in
sections 3 and 4. Finally, in Appendices A and B we display some
useful properties of orientable theories.

\noindent {\bf Acknowledgements.}
The first author is in debt to A.Merkurjev for inspiring
discussions at the initial stage of the work.
He is especially grateful to the Institute for Advanced Study (Princeton)
for excellent working conditions.

The main result of the paper was obtained during the stay of the
second author at Universit\"at Essen and the current text was
mostly written during his short-time visits to Universit\"at
Bielefeld and IH\'ES (Bures-sur-Yvette). The second author is very
grateful to all of these institutes for shown hospitality and
excellent working possibilities during these visits.

\noindent {\bf Notation.} Throughout the paper we use Greek
letters to denote elements of
cohomology groups and Latin for homological ones;
\begin{itemize}
\item $\Sm$ is a category of smooth quasi-projective algebraic varieties over a
field $k$.

\item $\Delta$ always denotes a diagonal morphism;

\item Symbol $\triv$ denotes trivial one-dimensional bundle;

\item For a vector bundle $\mathcal E$ over $X$ we write
$s(\mathcal E)$ for its section sheaf;

\item For a vector bundle $\mathcal E$ over $X$ we write
$\mathcal E^{\vee}$ for the dual to $\mathcal E$;

\item $\P(\mathcal E)=\text{Proj}(\text{Sym}^{*}(s(\mathcal E^{\vee})))$ is
the projective bundle of lines in $\mathcal E$;

\item $0\colon=\textrm{the point}\ [0\colon 0\colon\dots \colon 0\colon 1] \in \P^n $

\item typically $\P^n$ is regarded as a hyperplane in $\P^{n+1}$

\item $\P^{\infty}$ is a space defined in \cite{Vo}

\item $\pt=\spe k$;

\end{itemize}

\noindent
For the convenience of perception we usually move indexes up and down
oppositely  to the predefined positions of $*$ or $!$.

\begin{section}{Some products in (co)homology}
Consider a symmetric $T$-spectrum $\mathcal A$~\cite[p.505]{Ja}
endowed with a multiplication $\mu\colon\mathcal A\wedge\mathcal
A\to\mathcal A$ making $\mathcal A$ a symmetric commutative ring
$T$-spectrum. Then the spectrum $\mathcal A$ determines bigraded
cohomology and homology theories on the category of algebraic
varieties~\cite[p.595]{Vo}. A ring structure in cohomology is then
given by the cup-product satisfying the following commutativity
law. For $\alpha\in A^{p,q}$ and $\beta\in A^{p^\prime,q^\prime}$,
one has:
\begin{eqnarray}
\alpha \cup \beta = (-1)^{pp^\prime}\epsilon^{qq^\prime}(\beta \cup \alpha),
\end{eqnarray}
where $\epsilon\colon A^{*,*} \to A^{*,*}$ is the involution
described in Appendix B.

Suppose that $\mathcal A$ is endowed with an element $\gamma \in
A^{2,1}(\P^{\infty})$ satisfying the following two conditions:
\begin{itemize}
\item[(i)]
$\gamma|_{\P^0}=0 \in A^{2,1}(\P^0)$,
\item[(ii)]
$\gamma|_{\P^1}=\Sigma_T(1) \in A^{2,1}_{\{0\}}(\P^1)$
is the $T$-suspension of the unit $1 \in A^{0,0}(\pt)$,
\end{itemize}
then the pair $(\mathcal A, \gamma)$ is called {\it an oriented
symmetric commutative $T$-spectrum}. If $\mathcal A$ can be
endowed with an element $\gamma \in A^{2,1}(\P^{\infty})$
satisfying the conditions $(\rm i)$ and $(\rm ii)$ then $\mathcal
A$ is called {\it an orientable symmetric commutative
$T$-spectrum}. For such a $T$-spectrum $\epsilon=\id$ by Lemma
\ref{ChernAndEpsilon} and the commutativity law is reduced to
$\alpha \cup \beta = (-1)^{pp'}(\beta \cup \alpha)$. In this case
it is convenient to set $A^0=\oplus_{p,q}A^{2p,q}$,
$A^1=\oplus_{p,q}A^{2p-1,q}$, $A_0=\oplus_{p,q}A_{2p,q}$,
 and $A_1=\oplus_{p,q}A_{2p-1,q}$,
where $A^{*,*}$ (resp. $A_{*,*}$) are (co)homology theories
represented by the $T$-spectrum $\mathcal A$.
 The functors
$A^*=A^0 \oplus A^1\colon\Sm\to \ab2$ and $A_*=A_0 \oplus
A_1\colon\Sm\to \ab2$ are (co)homology theories taking values in
the category of $\mathbb Z/2$-graded abelian groups. Although all
our duality results hold for bigraded (co)homology groups, we
shall work, for simplicity, with the $\mathbb Z/2$-grading just
introduced.

Multiplicativity of the $T$-spectrum $\mathcal A$ provides a canonical way~\cite[13.50]{Sw} to
 supply the functors $A^*$ and $A_*$ (contravariant and covariant,
respectively) with a product structure consisting  of
two cross-products
\begin{eqnarray*}
\underline\times\colon A_p(X)\otimes A_q(Y)\to A_{p+q}(X\times Y),\hskip 5mm
\overline\times\colon A^p(X)\otimes A^q(Y)\to A^{p+q}(X\times Y)
\end{eqnarray*}
and two slant-products
\begin{eqnarray*}
/\colon A^p(X\times Y)\otimes A_q(Y) \to A^{p-q}(X),\hskip 5mm
\backslash\colon A^p(X)\otimes A_q(X\times Y) \to A_{q-p}(Y).
\end{eqnarray*}
One also defines two inner products
 \begin{eqnarray*}
\cup: A^p(X)\otimes A^q(X) \to A^{p+q}(X),\hskip 5mm
\cap: A^p(X)\otimes A_q(X) \to  A_{q-p}(X),
\end{eqnarray*}
as $\alpha\cup\beta=\Delta^*(\alpha\overline\times\beta)$ and
$\alpha\cap a=\alpha\backslash\Delta_*(a)$, correspondingly.
The cup-product makes the group $A^*(X)$ an associative
skew-commutative $\mathbb Z/2$-graded unitary ring and this
structure is functorial. (Skew-commutativity is not obvious and implied by the
 orientability of $\mathcal A$ as it is shown in Appendix~\ref{aux2}).
The cap-product makes the group $A_*(X)$ a unital $A^*(X)$-module
($1\cap a=a$ for every $a\in A_*(X)$) and this structure
is functorial in the sense that $\alpha\cap
f_*(a)=f_*(f^*(\alpha)\cap a)$.

We shall need below the following associativity relations, which
are completely analogous to ones existing in the topological
context (see, for example,~\cite[13.61]{Sw}).
 For $\alpha\in A^*(X\times Y)$, $\beta\in A^*(Y)$, $\gamma\in A^*(X)$,
$a\in A_*(Y)$,  and $b\in A_*(X)$, we have:
\def\theenumi{{\bf A.\arabic{enumi}}}
\begin{enumerate}
\item $\alpha/(\beta\cap a)=(\alpha\cup p_Y^*(\beta))/a$\label{as1}
\item $\gamma\cup(\alpha/a)=(p^*_X(\gamma)\cup\alpha)/a$\label{as2}
\item $(\alpha/a)\cap b=p^X_*((\alpha\cap(a\underline\times b))$,\label{as3}
\end{enumerate}
where $p_X$ and $p_Y$ denote the corresponding projections.

We shall also need the following functoriality property of the
$/$-product (comp.~\cite[13.52.{\it iii}]{Sw}). For morphisms
$f\colon X\to X^\prime$, $g\colon Y\to Y^\prime$, and elements
$\alpha\in A^*(X^\prime \times Y^\prime)$ and $a\in A_*(Y)$, one
has: $(f\times g)^*(\alpha)/a=f^*(\alpha/g_*(a))$.

For the final object $\pt$ in $\Sm$ one, clearly, has
$A^*(\pt)=A_*(\pt)$. This provides us with a distinguished element
$\lb \pt\rb\in A_0(\pt)$ (fundamental class of the point) such
that for any smooth $X$ and arbitrary $\alpha\in A^*(X)$, one has:
$\alpha/\lb\pt\rb=\alpha$. (Here we assume the standard
identification $X\times\pt=X$.) One can easily verify that the
canonical isomorphism $A^*(\pt)=A_{*}(\pt)$ may be written as
$\alpha\mapsto\alpha\cap\lb\pt\rb$. Throughout the paper we
implicitly use this construction and usually denote $\lb\pt\rb$ by
$1$.

\end{section}

\begin{section}{Poincar\'e Duality Theorem}
Let $\mathcal A$ be an orientable symmetric commutative ring
$T$-spectrum. Then the involution $\epsilon$ from $(1.1)$
coincides with the identity as explained in Appendix~\ref{aux2}.
So that, the commutativity law is reduced to $\alpha \cup \beta =
(-1)^{pp'}(\beta \cup \alpha)$. Setting
$A^0=\oplus_{p,q}A^{2p,q}$, $A^1=\oplus_{p,q}A^{2p-1,q}$, we see
that the functor $A^*:=A^0 \oplus A^1$ takes value in the category
of skew-commutative $\mathbb Z/2$-graded rings. For what follows
it is convenient to give the following
\begin{de}
Let $\mathcal A$ be an orientable symmetric commutative ring
$T$-spectrum. A Chern element is an element $\gamma \in
A^0(\P^{\infty})$ such that $\gamma|_{\P^0}=0$ and the family
$\{1, \gamma|_{\P^1}\} \subset A^0(\P^1)$ is a free basis of the
free rank two $A^0(\pt)$-module $A^0(\P^1)$. Another term, which
can be used for a Chern element, is a non-homogeneous orientation
of $\mathcal A$.
\end{de}
 A Chern element $\gamma$ lifts to a Chern structure in the
cohomology theory $A^*$ in the sense of ~\cite[Def.3.2]{Pa3} and
to a commutative Chern structure in the homology theory $A_*$
(resp.~\cite[Def's.2.1.1, 2.2.12]{Pi}). In fact, for every line
bundle $L$ over $X \in \Sm$ there exists a diagram of the form $X
\xleftarrow{p} X^{\prime} \xrightarrow{f} \P(V)$, where $V$ is a
finite dimensional $k$-vector space, $X^{\prime}$ is a torsor
under a vector bundle over $X$, and  the morphism $f$ is such that
the line bundles $p^{*}(L)$ and $f^{*}(\mathcal O_V(-1))$ are
isomorphic ~\cite[3.23]{Pa3}. Denote by $c(L)$ the class in
$A^0(X)$ such that $p^*(c(L))=f^*(\gamma|_{\P(V)})$. By \cite{So}
the element $c(L)$ is well-defined and the assignment $L \mapsto
c(L)$ is a Chern structure in $A^*$. Moreover, the family of
operators ${c(L)\cap} : A^Z_*(X) \to A^Z_*(X)$ forms a commutative
Chern structure in the homology theory $A_*$.

Any Chern structure in $A^*$ (resp. on $A_*$) lifts to a trace
structure in the cohomology (resp. homology),
see~\cite[Thm.4.1.2]{Pa2} (resp.~\cite[Thm.5.1.4]{Pi}). Namely, to
every projective morphism $f: Y \to X$ of smooth varieties over
$k$ one assigns two grade-preserving operators $f_!: A^*(Y) \to
A^*(X)$ and $f^!: A_*(X) \to A_*(Y)$ satisfying a list of natural
properties. Precise definitions of trace structures in a ring
(co)homology theory is given in ~\cite{Pa2, Pi}. The operators
$f_!$ and $f^!$ are called trace operators. (By historical reasons
they called {\it integrations} in~\cite{Pa2}.) The trace
structures $f \mapsto f_!$ and $f \mapsto f^!$ are explicit and
unique up to the following normalization condition. For a smooth
divisor $i\colon D \hookrightarrow X$:
\begin{equation}
\label{NormalizationCoh} i_{!} \circ i^{*}= i_!(1) \cup\colon
A^{*}(X) \to A^{*}(X),
\end{equation}
\begin{equation}
\label{NormalizationHom} i_{*} \circ i^{!}= i_!(1) \cap \colon
A_{*}(X) \to A_{*}(X),
\end{equation}
and $i_!(1)= c(\mathcal L(D))$.

For a projective morphism $f\colon Y \to X$ the map $f_!: A^{*}(Y)
\to A^{*}(X)$ is a two-side $A^{*}(X)$-module homomorphism, {\it
i.e.}
\begin{eqnarray}
\label{coproj} f_!(f^{*}(\alpha) \cup \beta)=\alpha \cup
f_!(\beta)\\ \nonumber f_!(\alpha \cup f^{*}(\beta))=f_!(\alpha)
\cup \beta.
\end{eqnarray}

\begin{de}
Let $\mathcal A$ be an orientable symmetric commutative ring
$T$-spectrum equipped with a Chern element $\gamma \in
A^0(\P^{\infty})$. For a smooth projective variety $X$ with the
structure morphism $\pi\colon X\to\pt$ we call $\pi^!(1)\in
A_0(X)$ the {\bf fundamental class} of $X$ in $A_*$ and denote it
by $\lb X\rb$.
\end{de}

\begin{rem}
Definitely, the class $[X]$ depends rather on the pair
$(A_*,\gamma)$ than on the $T$-spectrum $\mathcal A$. However, we
prefer to use this simplified notation, since we always keep in
mind one chosen and fixed Chern element $\gamma$ throughout the
paper.
\end{rem}

 With the notion of fundamental class in hands, one can define duality maps
\begin{equation}
\mathcal D^\bullet\colon A^*(X)\to A_{*}(X)\text{\ \ \ as\ \ \ }
\mathcal D^\bullet(\alpha)=\alpha\cap\lb X\rb
\end{equation}
and
\begin{equation}
\mathcal D_\bullet\colon A_*(X)\to A^{*}(X) \text{\ \ \  as\ \ \ } \mathcal
D_\bullet(a)=\Delta_!(1)/a.
\end{equation}

\begin{theorem}[\bf Poincar\'e Duality]
\label{PoincareDuality}
Let $\mathcal A$
be an orientable symmetric commutative ring $T$-spectrum
equipped with a Chern element
$\gamma \in A^0(\P^{\infty})$.
Then for every smooth projective variety $X$
the maps $\mathcal D^\bullet$ and $\mathcal D_\bullet$ are mutually inverse isomorphisms.
\end{theorem}

If $\gamma \in A^{2,1}$ and $X$ is equi-dimensional of dimension $d$
then $[X] \in A_{2d,d}(X)$. In this case the isomorphism
$\mathcal D^\bullet$ identifies $A^{p,q}$ with $A_{2d-p,d-q}$.
One can extract the following nice consequence of the
Poincar\'e Duality theorem, which enables us to interpret trace maps in a way topologists
like to.

\begin{cor}
For projective varieties $X,Y\in\Sm$ and a morphism $f\colon X\to
Y$, one has:
\begin{equation*}
f_!=\mathcal D_\bullet^Y f_* \mathcal D^\bullet_X\hskip 1cm\text{and}\hskip 1cm
f^!=\mathcal D^\bullet_X f^* \mathcal D_\bullet^Y,
\end{equation*}
where $\mathcal D_X$ and $\mathcal D_Y$ are introduced above duality operators
 for varieties $X$ and $Y$, respectively.
\end{cor}
\Proof
To proof the first equality, one should just check that $f_*\mathcal D^\bullet_X=\mathcal D^\bullet_Yf_!$.
Taking into account that $[X]=f^![Y]$, one immediately derives the desired relation from the
 First Projection Formula below (Theorem~\ref{ProjectionFormula}).
The second statement can be proven in a similar way, but requires the ``dual'' projection formula
that we do not consider here.
\Endproof
The proof of Theorem~\ref{PoincareDuality} is based on  two projection formulae for cap- and slant-products.

\begin{theorem}[\bf The First Projection Formula]
\label{ProjectionFormula} For $X,Y\in\Sm$, a projective morphism
$f: Y \to X$, and any elements $\alpha \in A^*(Y)$ and $a \in
A_*(X)$, the relation
\begin{equation}
\label{TheFirstProjection}
f_*(\alpha\cap f^!(a))=f_!(\alpha)\cap a
\end{equation}
holds in the group $A_*(X)$.
\end{theorem}

We need a few simple corollaries of this theorem.

\begin{cor}
\label{Transposition}
Let
$\tau: X \times X \to X \times X$
be the transposition morphism. Then for any elements
$\alpha \in A^*(X)$, $\beta\in A^*(X \times X)$,
and
$a \in A_*(X \times X)$,
we have:
\def\theenumi{{\bf{\alph{enumi})}}}
\begin{enumerate}
\item\label{m1a}
\begin{equation*}
\Delta_!(\alpha) \cap a=\Delta_!(\alpha)\cap \tau_*(a)\\
\end{equation*}
\item\label{m1b}
\begin{equation*}
\Delta_!(\alpha) \cup \beta = \Delta_!(\alpha) \cup \tau^*(\beta)
\end{equation*}
\end{enumerate}
in  $A_*(X \times X)$ ($A^*(X \times X)$, respectively).
\end{cor}

\Proof
Consider the Cartesian square
\begin{equation}
\xymatrix{
X\ar^-\Delta[r]\ar_\id[d]&X\times X\ar^\tau[d]\\
X\ar^-\Delta[r]&X\times X.}
\end{equation}
Since the map $\tau$ is flat, the square is transversal due
to~\cite[B.7.4.]{Fu}. By the base change
property~\ref{Basechange}, one has: $\Delta^! \circ \tau_*
=\Delta^!$. By Theorem~\ref{ProjectionFormula}, one has:
\begin{equation*}
\Delta_!(\alpha) \cap a = \Delta_*(\alpha \cap \Delta^!(a))=
\Delta_*(\alpha \cap \Delta^!(\tau_*(a)))=
\Delta_!(\alpha) \cap \tau_*(a)
\end{equation*}
that implies~\ref{m1a}. To get~\ref{m1b} one  uses cohomological projection
formula~(\ref{coproj}) instead.
\Endproof

\begin{theorem}[\bf The Second Projection Formula]
\label{SlantProductAndPushForward} Let $f: Y \to X$ be a
projective morphism of  smooth varieties. Let also $T\in\Sm$.
Then for every $\alpha \in A^*(T \times Y)$ and $a \in A_*(X)$,
one has (in $A^*(T)$):
\begin{equation}\label{SlantRelation}
\alpha / f^!(a) = F_!(\alpha) / a,
\end{equation}
where $F=\id\times f$.
\end{theorem}

\begin{cor}
\label{SlantProductWithDiagonal} Let $X$ be a smooth projective
variety. Then in $A^*(X)$, we have:
\begin{equation}
\label{}
\Delta_!(1) / [X] = 1.
\end{equation}
\end{cor}
\Proof Denote by $p\colon X\to\pt$ the structure morphism and let
$P=\id\times p\colon X\times X\to X$ be the projection. By
Theorem~\ref{SlantProductAndPushForward}, one has:
\begin{equation}
\Delta_!(1)/ [X] = \Delta_!(1)/p^!(1)=P_!(\Delta_!(1))=1.
\end{equation}
\Endproof

Now we derive the main result as an easy consequence of
Corollaries~\ref{SlantProductWithDiagonal} and~\ref{Transposition}.

\noindent {\it Proof of Theorem~\ref{PoincareDuality}.} Let
$p_1,p_2\colon X\times X\to X$ denote corresponding projections.
Observe that for every $\gamma \in A^*(X \times X)$ one has the
relation $\Delta_!(1) \cup \gamma = \gamma \cup \Delta_!(1)$. (In
fact, the element $\Delta_!(1)$ is of degree zero, because the map
$\Delta_!(1)$ is grade-preserving.) Thus, one has:
\begin{eqnarray}
\Delta_!(1)/(\alpha \cap\lb X\rb)\overset{\ref{as1}}=&
(\Delta_!(1)\cup p_2^*(\alpha))/\lb X\rb
\overset{\ref{Transposition}.\ref{m1b}}=&(\Delta_!(1) \cup
p_1^*(\alpha))/\lb X\rb \\\nonumber
=&(p_1^*(\alpha) \cup \Delta_!(1))/ [X]\overset{\ref{as2}}=&\alpha\cup(\Delta_!(1)/\lb X\rb)=\alpha.
\end{eqnarray}

On the other hand, using~\ref{Transposition}{.a}, one has:
\begin{eqnarray}
(\Delta_!(1)/ a )\cap\lb X\rb&\overset{\ref{as3}}=&
p_*(\Delta_!(1)\cap(a \underline\times \lb X\rb))
=p_*(\Delta_!(1)\cap(\lb X\rb \underline\times a))\nonumber\\
&\overset{\ref{as3}}=&(\Delta_!(1)/\lb X\rb)\cap a=a.
\end{eqnarray}
\hfill\qed
\bigskip

To complete the prove of Theorem \ref{PoincareDuality} one needs
to check formulae~\ref{TheFirstProjection}
and~\ref{SlantRelation}.
\end{section}

\begin{section}{Proof of The First Projection Formula.}

It is convenient to introduce
a class $\mathfrak V$ of projective morphisms
$f: Y \to X$ for which the relation
\begin{equation}
\label{KeyProjectionFormula}
f_*(\alpha \cap f^!(a))=f_!(\alpha) \cap a
\end{equation}
holds in $A_*(X)$ for every elements $\alpha \in A^*(Y)$ and $a
\in A_*(X)$.

Obviously, this class is closed with respect to composition.

We prove Theorem~\ref{ProjectionFormula} in several stages showing consequently that
the following classes of morphisms are contained in the class $\mathfrak V$.

\begin{itemize}
\item  Zero-section morphisms of  line bundles:
$s\colon Y\mono \P(\triv\oplus \L)$;
\item  Closed embeddings $i\colon D\mono X$ of smooth divisors;
\item Zero-sections of a finite sum of line bundles:
\begin{equation*}
s\colon Y\mono \P(\triv\oplus \L_1\oplus \L_2\oplus\cdots\oplus \L_n));
\end{equation*}
\item Zero-sections of arbitrary vector bundles:
$s\colon Y\mono \P(\triv\oplus \mathcal V)$;
\item Closed embeddings $i\colon Y\mono X$;
\item Projections $p\colon X\times\P^n\to X$;
\end{itemize}

\begin{lemma}\label{onedim}
Let $\L$ be a line bundle over a smooth variety Y.
Then the zero-section
$s\colon Y \mono \P(\triv\oplus \L)$
belongs to $\mathfrak V$.
\end{lemma}

\Proof
 The map $s$ is a section of the
projection map $p\colon \P(\triv\oplus \L)\to Y$.
Let $\alpha\in A^*(Y)$ and $a\in
A_*(\P(\triv\oplus \L))$. The desired relation follows from
(\ref{NormalizationHom}) and (\ref{NormalizationCoh}):
\begin{eqnarray}
s_*(\alpha\cap s^!(a))&=s_*(s^*p^*(\alpha)\cap s^!(a))&=
p^*(\alpha)\cap s_*s^!(a)\\\nonumber
&=p^*(\alpha)\cap\left(s_!(1)\cap
a\right)&=s_!(s^*p^*(\alpha)) \cap a = s_!(\alpha) \cap a.
\end{eqnarray}

\Endproof

\begin{prop}
\label{reduction}
Let $X,Y\in\Sm$,
$i\colon Y\mono X$
be a closed embedding with a normal bundle $\N$.
If the zero-section morphism
$s\colon Y\mono \P(\triv\oplus\N)$
belongs to $\mathfrak V$ then $i$
belongs to $\mathfrak V$.
\end{prop}

\Proof
Consider the following deformation diagram,
in which $B$ is the blowup of $X\times\A^1$ at
$Y\times\{0\}$. This diagram has transversal squares.

\begin{equation}
\label{KeyDiagram} \xymatrix{
&B-Y\times\A^1\ar[d]_{k_B }\\
\P(\triv\oplus\N)\ar@{^(->}[r]_{k_0}\ar@/^0.7pc/[d]^p
&B&X\ar@{_(->}[l]^{k_1}\\
Y\ar@{^(->}[r]_{j_0}\ar@{^(->}[u]^{s}
&Y\times\A^1\ar@{^(->}[u]^{t}
&Y\ar@{_(->}[l]^{j_1}\ar@{^(->}[u]^{i}}
\end{equation}

One can easily see that the left-hand part of our diagram satisfies the conditions
of Lemma~\ref{useful}.

First, we shall show that the morphism $t$ in
Diagram~\ref{KeyDiagram} belongs to the class $\mathfrak V$.
Let $\alpha \in A^*(Y \times \A^1)$ and $a \in A_*(B)$.
Using Lemma~\ref{useful} we can rewrite $a$ as
$k^B_*(a_B)+k^0_*(a_0)$, where $a_0\in A_*(\P(\triv\oplus\N))$
and  $a_B\in A_*(B-Y\times\A^1)$. From the Gysin exact sequence,
we have:
\begin{eqnarray}
t^!k^B_*&=&0\quad\text{and}\\\label{2.5}k^*_Bt_!&=&0.
\end{eqnarray}
Therefore, $t_*(\alpha\cap t^!k^B_*(a_B))=0$ and $t_!(\alpha)\cap
k^B_*(a_B)=0$. (The second relation yields
from~\ref{2.5}:
$t_!(\alpha)\cap k^B_*(a_B)=k^B_*(k^*_Bt_!(\alpha) \cap a)=0$.)
Thus, one has:
\begin{equation}
\label{IntermidiateFinal} t_*(\alpha \cap t^!(a)) =t_*(\alpha \cap
t^!k^0_*(a_0)).
\end{equation}
Applying Lemma~\ref{lifting} to the left-hand-side square of
Diagram~\ref{KeyDiagram} and denoting $j_0^*(\alpha)$ by
$\alpha_0$, one has:
\begin{equation}
\label{IntermidiateOne} t_*(\alpha\cap
t^!k^0_*(a_0))=k^0_*s_*\left(\alpha_0\cap s^!(a_0)\right).
\end{equation}
Similarly:
\begin{equation}
\label{IntermidiateTwo} t_!(\alpha)\cap
a=k^0_*\left(s_!(\alpha_0)\cap a_0\right).
\end{equation}
By the proposition assumption, we have the relation
$s_*(\alpha_0\cap s^!(a_0))=s_!(\alpha_0)\cap a_0$. Combining this
with equalities \ref{IntermidiateFinal}, \ref{IntermidiateOne},
and \ref{IntermidiateTwo}, one gets:
\begin{equation}
\label{IntermidiateRelation}
t_*(\alpha \cap t^!(a)) = t_!(\alpha) \cap a.
\end{equation}
We now move the desired relation one more step further to the
right in Diagram~\ref{KeyDiagram} and show that $i\in\mathfrak V$.
Observe that  $k^1_*$ is a monomorphism. Therefore, it suffices to
check that
 for every elements $\alpha_1 \in A^*(Y)$ and $a_1 \in A_*(X)$ we have:
\begin{equation}
k^1_* i_*(\alpha_1\cap i^!(a_1))=k^1_*(i_!(\alpha_1)\cap a_1).
\end{equation}
Setting $\alpha=(j^*_1)^{-1}(\alpha_1)\in A^*(Y \times \A^1)$,
$a=k^1_*(a_1) \in A_*(B)$, and applying Lemma~\ref{lifting} to the
right-hand-side square of Diagram~\ref{KeyDiagram}, one has:
$k^1_*i_*(\alpha_1\cap i^!(a_1))=t_*(\alpha\cap t^!(a))$. In the
same way: $k^1_*(i_!(\alpha_1)\cap a_1)= t_!(\alpha) \cap
k^0_*(a_0)= t_!(\alpha)\cap a$. Combining these two relations with
\ref{IntermidiateRelation}, one sees that $i\in\mathfrak V$.
\Endproof
\begin{cor}
\label{div} For a smooth divisor $i\colon D\mono X$ the morphism
$i$ lies in $\mathfrak V$.
\end{cor}

\begin{cor}
\label{splbdl}
Let
$\mathcal W=\L_1\oplus\cdots\oplus \L_n$
be an $n$-dimensional vector bundle over a variety
$Y$ which splits in the sum of line bundles.
Then the zero-section morphism
$s\colon Y\mono\P(\triv\oplus\mathcal W)$
 belongs to the class $\mathfrak V$.
\end{cor}

\Proof
Apply Corollary~\ref{div} to each step of the filtration
\begin{equation}
Y\overset{i_1}\mono\P(\triv\oplus \L_1)\overset{i_2}\mono
\cdots\overset{i_n}\mono\P(\triv\oplus\mathcal W),
\end{equation}
where the morphisms $i_j$ are zero-sections of $\L_j$.
\Endproof

In order to proceed with the case of an arbitrary vector-bundle,
we need the
homological analogue of the splitting principle.
Consider
a vector bundle $\mathcal E\to Y$  of constant rank $n$
over a smooth irreducible variety $Y$. Let
$\GL_n$
be the corresponding principal
$GL_n$-bundle over $Y$,
$T_n \subset GL_n$
be the diagonal tori, and
$Y^{\prime}=\GL_n / T_n$
be the orbit variety with the projection morphism
$p: Y^{\prime} \to Y$.
Finally, we denote by
$\E^{\prime}=\E \times_Y Y^{\prime}$
 the pull-back of the vector bundle $\E$.

\begin{prop}
\label{splitting} The bundle $\E^\prime$ splits in a direct sum of
line bundles and the map $p_*\colon A_*( Y^\prime)\to A_*(Y)$ is a
universal splitting epimorphism (i.e. for any base-change
$Z\to Y$ the induced map $A_*(Z\times_Y Y^\prime)\to A_*(Z)$ is a
splitting epimorphism).
\end{prop}
\Proof
The projection
$\mathcal G\mathcal L_n \to Y^{\prime}$
and the natural
$T_n$-action on
$\GL_n$
makes it a principal $T_n$-bundle over
$Y^{\prime}$.
Moreover, if
$\GL^{\prime}_n=
\GL_n \times_Y Y^{\prime}$
is the pull-back of
$\GL_n$, there is a natural isomorphism of principal
$GL_n$-bundles
\begin{equation}
\GL_n \times_{T_n} GL_n \to
\GL_n^{\prime}
\end{equation}
over $Y^{\prime}$. The bundle
$\E^{\prime}$ over $Y^{\prime}$
corresponds exactly to the principal
$GL_n$-bundle
$\GL_n^{\prime}$.
Thus, the mentioned isomorphism of principal
$GL_n$-bundles over $Y^{\prime}$
shows that the bundle $\E^{\prime}$
splits in a direct sum of line bundles
(say corresponding to the fundamental
characters $\chi_1, \chi_2, \dots, \chi_n$ of the tori $T_n$).
This proves the first assertion of the proposition.

To prove the second one, consider a Borel subgroup $B_n$ in $GL_n$
(say the subgroup of all upper triangle matrices) and let $U_n$ be
the maximal unipotent subgroup of $B_n$ (the group of upper
triangle matrices with 1's on the diagonal). Let $\mathcal F=\GL_n
/ B_n$ (this is just the flag bundle over $Y$ associated to $\E$).
The bundle $\mathcal F$ comes equipped with projections $q:
\mathcal F \to Y$ and $r: Y^{\prime} \to \mathcal F$, where the
projection $r$ is induced by the inclusion $T_n \subset B_n$.
Using the natural $U_n$-action on $\GL_n$,  it is easy to check
that there is a tower of morphisms:
\begin{equation}
\GL_n =
S_m \to S_{m-1} \to \dots \to S_1 =
\mathcal F,
\end{equation}
which has a principal $\mathbb G_a$-bundle on each level (each
level is a torsor over the trivial rank one vector bundle). By the
strong homotopy invariance property~\cite[2.2.6]{Pa3}, the induced
map on homology $r_*\colon A_*(Y^{\prime}) \to A_*(\mathcal F)$ is
an isomorphism.

As it was already mentioned,
$\mathcal F$ is a full flag bundle over $Y$  associated to the
bundle $\E$. Thus, there is a tower of morphisms
\begin{equation}
\mathcal F = Z_s \to Z_{s-1} \to \dots \to Z_1 = Y
\end{equation}
in which each level is a projective bundle
associated to a vector bundle. By the
Projective Bundle Theorem (PBT)~\ref{PBT2}, we have a split epimorphism
in homology induced on each floor.
Therefore, the map $q_*\colon A_*(\mathcal F) \to A_*(Y)$ is a split epimorphism as well.

These proves that the map
$p_*: A_*(Y^{\prime}) \to A_*(Y)$
is also an epimorphism.

One can easily check that all necessary properties of the morphisms
$p,q$, and $r$ are base-change invariant. Therefore, the constructed
splitting epimorphism is universal.
\Endproof

\begin{prop}
\label{vectorbdls}
Let
$s\colon Y\mono\P(\triv\oplus \mathcal V)$
be  the zero-section of
the finite-dimensional vector bundle
$\mathcal V$.
Then $s\in\mathfrak V$.
\end{prop}

\Proof
Letting $Y^{\prime}$ be as above, denote by $\mathcal V^{\prime}$
the pull-back of the bundle $\mathcal V$
with respect to the morphism $p$.
Then by Proposition~\ref{splitting} the bundle $\mathcal V^{\prime}$
splits in a direct sum of line bundles and the induced map
\begin{equation}
\bar p_*\colon A_*(\P(\triv\oplus\mathcal V^{\prime})) \to
A_*(\P(\triv\oplus\mathcal V))
\end{equation}
is an epimorphism.

 Let $s: Y \to
\P(\triv\oplus\mathcal V)$ and $\bar s\colon Y^{\prime} \to
\P(\triv\oplus\mathcal V^{\prime})$ be morphisms induced by
zero-sections of the corresponding vector bundles. Then the diagram
\begin{equation}
\xymatrix{\P(\triv\oplus\mathcal V^{\prime})\ar[r]^{\bar p}
&\P(\triv\oplus\mathcal V)\\
Y^{\prime}\ar^{\bar s}[u]\ar^{p}[r]&Y\ar^{s}[u]}
\end{equation}
is transversal.

Let $\alpha\in A^*(Y)$ and $a\in A_*(\P(\triv\oplus\mathcal
V))$. Choosing $b\in A_*(\P(\triv\oplus\mathcal V^\prime))$
such that $a=\bar p_*(b)$ and applying Lemma~\ref{lifting}, one gets:
\begin{equation}
s_*(\alpha\cap s^!(a))= {\bar
p}_*\bar s_*(p^*(\alpha)\cap \bar s^!(b))
\end{equation}
and
\begin{equation}
s_!(\alpha)\cap a={\bar p}_*(\bar s_!p^*(\alpha)\cap b).
\end{equation}
Two expressions on the right-hand-sides coincide by Proposition~\ref{splbdl}.
\Endproof

\begin{cor}
\label{TheCaseOfClosedImbedding}
Let
$i\colon Y\mono X$ be a closed embedding. Then,
$i\in\mathfrak V$.
\end{cor}

\Proof
Applying Proposition
\ref{reduction}
we reduce the question to the case of the zero-section morphism
$s\colon Y\mono \P(\triv\oplus\N)$
 of the normal bundle
$\N=\N_{X/Y}$.  The morphism $s$ belongs to $\mathfrak V$ by
Proposition~\ref{vectorbdls}.
\Endproof

In order to check that  projection morphisms $p: X \times \P^* \to
X$ belong to $\mathfrak V$ we need a few auxiliary results
(\ref{diamond}--\ref{UpThenDownInHomology}).

\begin{nota}
For a projective morphism $f$ we denote, from now on, the map
$f_*f^!$ by $f^\dia$ and $f_!f^*$ by $f_\dia$.
\end{nota}

\begin{lemma}\label{diamond}
\def\theenumi{{\bf \alph{enumi})}}
\begin{enumerate}
\item\label{diamond1}(left distributivity) Let $a,b,c$, and $p$ be projective morphisms.  If $a^\dia=b^\dia+c^\dia$ then
$(pa)^\dia=(pb)^\dia+(pc)^\dia$, provided that both sides of the equality are well defined.
\item\label{diamond2} Given a transversal square with projective morphisms
 $f$ and $g$
$$\xymatrix{
X\times_Z Y\ar[r]\ar_F[d]\ar@{..>}^h[dr]&Y\ar^f[d]\\
X\ar_g[r]&Z }
$$
one has the following equalities: $h^\dia=g^\dia f^\dia=f^\dia g^\dia$.
\item\label{diamond3}
In the square above: $g_*F^\dia=f^\dia g_*$.
\item \label{diamond4}
Let $s_i$ be the standard embedding $\P^{n-i}\hookrightarrow\P^n$ and $p_n\colon\P^n_X\to X$
be the projection map.
Let $\psi_i$ be the same as in the Projective Bundle Theorem
(see~\ref{PBT2}).
Then $p_*^ns_i^\dia=\psi_i$.
\end{enumerate}
\end{lemma}
\Proof
Part~\ref{diamond1} immediately follows from the definition of the operation $\dia$,
~\ref{diamond2} and~\ref{diamond3} are trivial corollaries of the transversal base-change property,
 ~\ref{diamond4} easily follows from the PBT.
\Endproof

Fix now a variety $X\in\Sm$ and take the $n$-dimensional
projective space $\P^n_X$  over $X$. (Up to the end
of this section all the schemes are considered over the base scheme $X$ and the
 product is implicitly taken over $X$.)
Due to the PBT, the element $\Delta_!(1)\in A^*(\P^n\times\P^n)$ may be decomposed as
\begin{equation}
\Delta_!(1)= 1\boxtimes\zeta^n+\zeta^n\boxtimes 1 + \sum_{i,j=1}^n
a_{ij}\zeta^i\boxtimes\zeta^j,
\end{equation}
where $\zeta=e(\mathcal O(1))$ is the canonical generator of
$A^*(\P^n)$ as an $A^*(X)$-algebra and $a_{ij}\in A^*(X)$ (see
~\cite[Lemma 1.9.3]{Pa}).

This equality together with the previous lemma gives us the
following decomposition of the identity operator $\id_{\P^n}$.
Taking into account the relation
$s_{ij}^\dia(x)=(\zeta^i\boxtimes\zeta^j)\cap x$, where
$s_{ij}\colon\P^{n-i}\times\P^{n-j}\hookrightarrow
\P^{n}\times\P^{n}$ is the standard embedding, we can rewrite the
cap-product with $\Delta_!(1)$ operator in the form:
\begin{equation}
\Delta^\dia=(\Delta_!(1)\cap)=s_{0n}^\dia+s_{n0}^\dia+\sum_{i,j=1}^n
a_{ij}s_{ij}^\dia.
\end{equation}
Applying Lemma~\ref{diamond} to transversal squares:
\begin{equation}
\xymatrix{\P^{n-i}\times\P^{n-j}\ar[d]\ar@{^(->}[r]\ar@{..>}_{p_{1,n}s_{ij}}[dr]&\P^n\times\P^{n-j}\ar^{p_{1,n-j}}[d]\\
\P^{n-i}\ar@{^(->}_{s_i}[r]&\P^n}
\end{equation}
(where we denote by $p_{1,k}$ the projection map $\P^n\times\P^k\to\P^n$),
one gets the following equality:
\begin{equation}\label{dec}
\id=(p_{1,n}\Delta)^\dia=\sum_{i,j=0}^n
a_{ij}(p_{1,n}s_{ij})^\dia=\sum_{i,j=0}^n
a_{ij}p_{1,n-j}^\dia s_{i}^\dia.
\end{equation}(here we set $a_{i0}=a_{0i}=0$, unless $i=n$.)

\begin{lemma}
\label{diag} For the projection morphism $p_n\colon \P^n\to X$,
we have:
\def\theenumi{{{\alph{enumi})}}}
\begin{enumerate}
\item\label{proj-a} $p_n^\dia=-\sum_{j=1}^n a_{nj}p_{n-j}^\dia;$
\item\label{proj-b} $p^n_\dia=-\sum_{j=1}^n a_{nj}p^{n-j}_\dia.$
\end{enumerate}
\end{lemma}

\Proof
Let us check the first statement. For $n=0$ we, trivially, have $p_0^\dia=\id$.
Applying the map $p^n_*$ to both sides of~\ref{dec} and then using Lemma~\ref{diamond}.c
for the transversal squares
\begin{equation}
\xymatrix{
\P^n\times\P^{n-j}\ar[r]\ar_{p_{1,n-j}}[d]&\P^{n-j}\ar^{p_{n-j}}[d]\\
\P^n\ar_{p_n}[r]&X,}
\end{equation}
one gets:
\begin{equation}
p_*^n=\sum_{i,j=0}^n
a_{ij}p_*^np_{1,n-j}^\dia s_{i}^\dia=\sum_{i,j=0}^n
a_{ij}p_{n-j}^\dia(p_n^*s_{i}^\dia).
\end{equation}
Taking into account that due to~\ref{diamond}.d $p^n_*s_i^\dia=\psi_i$,
($p^n_*=\psi_0$), one has:
\begin{equation}
\label{FinalRel}
\left(p_n^\dia+\sum_{j=1}^n a_{nj} p_{n-j}^\dia\right)\psi_n =
- \sum_{i,j=1}^{n-1}(\cdots)\psi_i.
\end{equation}
By the  PBT, for any $x\in A_*(X)$ we can choose an element
$\phi(x)\in A_*(\P^n_X)$ such that $\psi_i(\phi(x))=\left\{
\begin{array}{ll}0,&i<n\\ x, & i=n.
\end{array} \right.$

Applying operators on both sides of~(\ref{FinalRel}) to $\phi(x)$,
we get:
\begin{equation}
0= p_n^\dia+\sum_{j=1}^n a_{nj} p_{n-j}^\dia
\end{equation}

This finishes the proof of case~\ref{proj-a}.
The cohomological relation~\ref{proj-b} may be proved by dualization of these
arguments or found in~\cite[Section 1.10]{Pa}.
\Endproof

\begin{prop}
\label{UpThenDownInHomology}
Let $p_n$ denote, as before, the projection morphism $p_n\colon \P^n_X\to X$. Then for every element $a \in A_*(X)$,
one has:
\begin{equation*}
p^n_*(p^!_n(a)) = p^n_!(1) \cap a.
\end{equation*}
\end{prop}

\Proof
Rewriting the Proposition statement in our notation,
  we should verify the relation
$p_n^\dia(a) = p^n_\dia(1) \cap a$.
We proceed by induction on $n$.
The case $n=0$ is trivial.
Let the proposition hold for $n<N$.
Then for $p_N$, by Lemma~\ref{diag}, we have:
\begin{equation}
p_N^\dia(a)=-\sum_{j=1}^N a_{Nj} p_{N-j}^\dia(a)
\end{equation}
and
\begin{equation}
p^N_\dia(1)\cap a=-\sum_{j=1}^N a_{Nj} p^{N-j}_\dia(1)\cap a
\end{equation}
By the induction hypothesis the expressions on the right-hand-side coincide.
The induction runs.
\Endproof

\begin{prop}
\label{TheCaseOfProjection}
For every integer $n\ge 0$ the projection morphism
$p\colon \P^n_X\to X$
belongs to the class
$\mathfrak V$.
\end{prop}

\Proof
Given $\alpha \in A^*(\P^n_X)$ and $a \in A_*(X)$ one
should verify that
\begin{equation}\label{ProjectionCase}
p_*(\alpha \cap p^!(a)) = p_!(\alpha) \cap a.
\end{equation}
Clearly, both sides of~(\ref{ProjectionCase}) are $A^*(X)$-linear.
By the PBT, $A^*(\P^n_X)$ is generated as an $A^*(X)$-module by
the elements $\zeta^j$. Thus, it suffices to check the Proposition
just for these elements. From \cite[Lemma 1.9.1]{Pa}, we have a
relation $\zeta^j=i^j_!(1)$ in $A^*(\P^n_X)$, where $i^j\colon
\P^{n-j}_X\mono \P^n_X$ is the standard embedding map and the
element $\zeta^j \in A^*(\P^n)$ is considered here as lying in
$A^*(\P^n_X)$ via the pull-back operator for the projection
$\P^n_X \to \P^n$. Denote by $p_j$ the projection map $\P^{n-j}_X
\to X$. Since $p\circ i^j=p_j$, we have by
Proposition~\ref{TheCaseOfClosedImbedding}:
\begin{equation}
p_*(\zeta^j\cap p^!(a))= p_*i^j_*(1\cap
i_j^!p^!(a))= p^j_*p_j^!(a).
\end{equation}
Using Proposition
\ref{UpThenDownInHomology}, one has:
\begin{equation}
p^j_*p^!_j(a)=p^j_!(1)\cap a=p_!i^j_!(1) \cap a=p_!(\zeta^j)\cap a.
\end{equation}
\Endproof

This finishes proof of Theorem~\ref{ProjectionFormula}.

\end{section}

\begin{section}{Proof of The Second Projection Formula}
The strategy of the proof of  Theorem~\ref{SlantProductAndPushForward}
is very similar to one used in the previous section.
It is again convenient to  introduce a class $\mathfrak W$
consisting of projective morphisms $f: Y \to X$ such that for any $T\in\Sm$, $\alpha\in A^*(T\times Y)$, and
$a\in A_*(X)$ the relation
\begin{equation}
F_!(\alpha)/a=\alpha/f^!(a)
\end{equation}
holds in $A^*(T)$. (Here $F=\id\times f$. Below we use similar notation rules.)

We show that the following classes of morphisms lie in $\mathfrak W$.
\begin{itemize}
\item Zero-sections  of vector bundles:
$s\colon Y\mono \P(\triv\oplus \V)$;
\item Closed embeddings $i\colon Y\mono X$;
\item Projections $p\colon X\times\P^n\to X$;
\end{itemize}
Since the class $\mathfrak W$ is closed with respect to composition, this will imply our
formula for all projective morphisms.
\begin{lemma}
\label{SlantForTheSectionOfProjectiveBundle}
Let $\V$ be a vector bundle over a smooth variety $Y$
and let
$s\colon Y\mono \P(\triv\oplus \V)$
be the zero-section of the projection $p\colon \P(\triv\oplus \V)\to Y$.
Then the morphism $s$ belongs to the class $\mathfrak W$.
\end{lemma}

\Proof
Let $\alpha \in A^*(T \times Y)$ and $a \in A_*(\P(\triv\oplus \V))$.
Using functoriality of
the slant-product, an associativity relation, and formula~\ref{coproj},
one gets:
\begin{eqnarray}
\alpha/s^!(a)&=&\alpha/p_*(s_!(1)\cap a)=
P^*(\alpha)/(s_!(1)\cap a)\nonumber\\
&\overset{\ref{as1}}=&(P^*\left(\alpha)\cup (1\times
s_!(1))\right)/a= \left(P^*(\alpha) \cup S_!(1)\right)/a=
S_!(\alpha)/a.
\end{eqnarray}
(Here the relation $1\times s_!(1)=S_!(1)$ appears from the
base-change property applied to the product with $T$.)
\Endproof

\begin{prop}
\label{SlantReduction}
Any closed embedding morphism $i :Y \hookrightarrow X$ of smooth varieties belongs to the class $\mathfrak W$.
\end{prop}

\Proof
Denote by $\P(\triv\oplus\N)$ the  projectivization corresponding to the normal
bundle $\N=\N_{X/Y}$. It is endowed with the zero-section morphism $s\colon Y\mono \P(\triv\oplus \N)$.

As well as in the proof of Theorem~\ref{ProjectionFormula} our
arguments are based on the deformation diagram which obtained
from~(\ref{KeyDiagram}) by multiplication with a variety $T\in\Sm$.
For convenience, we reproduce this diagram here.

\begin{equation}
\xymatrix{
&T \times B-T \times Y\times\A^1\ar[d]_{K_B }\\
T \times \P(\triv\oplus\N)\ar@{^(->}[r]_{K_0}&T \times B&T\times X\ar@{_(->}[l]^{K_1}\\
T \times Y\ar@{^(->}[r]_{J_0}\ar@{^(->}[u]^{S}&T \times Y\times\A^1\ar@{^(->}[u]^{I_t}&
T \times Y\ar@{_(->}[l]^{J_1}\ar@{^(->}[u]^{I}}
\end{equation}

First of all, we show that $I_t\in\mathfrak W$. Namely, we should
prove  that for any elements $\alpha\in A^*(T\times Y\times\A^1)$
and  $a\in A_*(B)$ the relation
\begin{equation}
\label{SlantIntermidiateRelation}
\alpha/i^!_t(a) = I^t_!(\alpha)/a.
\end{equation}
holds in $A^*(T)$.

Exactly as in the proof of Theorem~\ref{ProjectionFormula} one can
 rewrite $a$ as a sum
$k^B_*(a_B)+k^0_*(a_0)$, where $a_0\in A_*(\P(\triv\oplus\N))$ and
$a_B\in A_*(B-Y\times\A^1)$ and obtain the equalities:
\begin{equation}
\label{SlantIntermidiateOne} \alpha / i_t^!(a)=\alpha /
i_t^!k^0_*(a_0)=\alpha_0 / s^!(a_0),
\end{equation}
where $\alpha_0=J_0^*(\alpha)$.

Similarly, one gets the relation:
\begin{equation}
\label{SlantIntermidiateTwo}
I^t_!(\alpha) / a=
S_!J^*_0(\alpha) / a_0 =
S_!(\alpha_0) / a_0.
\end{equation}
By Lemma~\ref{SlantForTheSectionOfProjectiveBundle} $\alpha_0 / s^!(a_0)=S_!(\alpha_0) / a_0$,
which proves~(\ref{SlantIntermidiateRelation}).

Since the map $J_1^*$ is an isomorphism, we can set
$\alpha=(J^*_1)^{-1}(\alpha_1)\in A^*(T \times Y \times \A^1)$ and
$a=k^1_*(a_1) \in A_*(B)$. Applying Lemma~\ref{lifting} again, one
gets:
\begin{eqnarray}
\label{SlantRelOne}
\alpha_1 / i^!(a_1)&=&J_1^*(\alpha) / i^!(a_1)=\alpha / i_t^!(a)\quad\text{and}\\
\label{SlantRelTwo} I_!(\alpha_1) / a_1&=&I_!J_1^*(\alpha) /
a_1=I^t_!(\alpha) / a.
\end{eqnarray}
Combining these equalities with relation~(\ref{SlantIntermidiateRelation})
 proves the proposition.
\Endproof

\begin{prop}
\label{proj1}
Let $X,T\in\Sm$,
 $p\colon X\times \P^n \to X$
be the projection morphism, and $P= \id \times p\colon T\times
X\times\P^n\to T \times X$. Then for every elements $\alpha\in
A^*(T\times X\times \P^n)$ and $a \in A_*(X)$, one has a relation:
\begin{equation}
\label{proj1Relation}
\alpha / p^!(a) = P_!(\alpha) / a
\end{equation} in $A^*(T)$.
\end{prop}

\Proof
Consider the following commutative diagram with transversal square:
\begin{equation}
\xymatrix{X\times\P^n\ar_{p=p_0}[dr]&X\times\P^{n-r}\ar_i[l]\ar^{p_r}[d]&T\times X\times \P^{n-r}\ar_{\bar q}[l]\ar^{P_r}[d]\\
&X&T\times X\ar_q[l]
}
\end{equation}
Clearly, both sides of~(\ref{proj1Relation}) are $A^*(T)$-linear.
So, we may assume that $\alpha=\zeta^r_{T\times X}$.
Since $\zeta_{T\times X}^r=I_!(1_{T\times X}) \in A^*(T \times
X\times\P^n)$, one has:
\begin{eqnarray}
\zeta_{T\times X}^r/p^!(a)&=&I_!(1_{T\times X})/ p^!(a)=1_{T\times X} / i^!p^!(a)\nonumber\\
&=&1/ p_r^!(a)=P^*_r(1)/p^!_r(a)=1/p^r_*p^!_r(a).
\end{eqnarray}
By Proposition~\ref{UpThenDownInHomology}:
\begin{equation}
1/p^r_*p^!_r(a)=1/(p^r_!(1_X)\cap a)=q^*p^r_!(1)/a.
\end{equation}
Applying the base-change property to the square in the diagram above, we get:
\begin{equation}
q^*p^r_!(1_X)=P^r_!(1_{T\times
X})=P_!(I_!(1))=P_!(\zeta^r_{T\times X}).
\end{equation}
The proposition is proven.
\Endproof
\end{section}
\appendix
\begin{section}{Some properties of a trace structure}
\label{aux} In this Appendix we give a brief description of some
useful properties of a trace structure, which are utilized in the
paper. Although we make a deal with both cohomological and
homological contexts, we present all the results for homology.
Cohomological variant is ``dual'' in the obvious sense and may be
found in~\cite{Pa}. All the proofs for the homological case not
provided below can be found in~\cite{Pi}.

We, first, give a definition of a transversal square
following A.Mer\-kurjev~\cite{Me}.
\begin{de}\label{transversal}
We call a square
\begin{equation*}
\xymatrix{Y^\prime\ar^{\bar f}[r]\ar_{\bar g}[d] & X^\prime\ar[d]^{g}\\
Y\ar_{f}[r]&X}
\end{equation*}
in the category $\Sm$ {\bf transversal} if

\def\theenumi{{\bf{\alph{enumi})}}}
\begin{enumerate}
   \item it is Cartesian in the category ${\bf Sch/}k$ of all schemes over the field $k$;

    \item \label{trans2} the following sequence of tangent bundles over $Y^\prime$ is exact:
\end{enumerate}
\begin{equation*}
0\to T_{Y^\prime}\overset{d\bar g\oplus d\bar f}{\to}{\bar g}^*T_Y\oplus {\bar f}^*T_{X^\prime}
\overset{dg-df}{\to} {\bar g}^*f^*T_X\to 0.
\end{equation*}
\end{de}

It is not hard to check that this definition is accordant to ones given in
~\cite[1.1.2]{Pa} or~\cite[1.1]{PY}.
Let us check, for example, that for a closed embedding $f$ condition
~\ref{trans2} implies
the isomorphism: $\bar g^*\mathcal N_{X/Y}\simeq \mathcal N_{X^\prime/Y^\prime}$.
The short exact sequence above may be viewed as a total complex of the bicomplex:
\begin{equation}
\xymatrix{0\ar[r]&{\bar g}^*T_Y\ar^{-df}[r]&{\bar g}^*f^*T_X\\
0\ar[r]&T_{Y^\prime}\ar_{d\bar f}[r]\ar^{d\bar g}[u]&{\bar f}^*T_{X^\prime}.\ar_{dg}[u]}
\end{equation}
Since~\ref{trans2} is exact, the bicomplex is acyclic. On the other hand, it is quasiisomorphic
to the two-term complex
$\bar g^*\mathcal N_{X/Y}\leftarrow \mathcal N_{X^\prime/Y^\prime}$.

\begin{property}[\bf Base-change for transversal squares]
\label{Basechange}
For any transversal square as above
with projective morphism $f$ the diagram
\begin{equation*}
\xymatrix{A_*(Y^\prime)\ar_{\bar g_*}[d]&A_*(X^\prime)\ar_{\bar
f^!}[l]\ar_{g_*}[d]
\\
A_*(Y)&A_*(X)\ar_{f^!}[l]}
\end{equation*}
commutes.
\end{property}

\begin{cor}
\label{lifting} Suppose, we are given a transversal square
\begin{equation*}
\xymatrix{X^\prime\ar^{g}[r] & X\\
Y^\prime\ar_{\bar g}[r]\ar^{\bar f}[u]&Y\ar[u]_f}
\end{equation*}
with projective morphism $f$. Let $\alpha\in A^*(Y)$ and $a\in
A_*(X^\prime)$. Then the following relations hold:
\def\theenumi{{\it \roman{enumi})}}
\begin{enumerate}
\item $f_*(\alpha\cap f^! g_*(a))= g_*\bar f_*(\bar g^*(\alpha)\cap\bar f^!(a))$
\item $f_!(\alpha)\cap g_*(a)= g_*(\bar f_!\bar g^*(\alpha)\cap a)$

\noindent Moreover, for a variety $T\in\Sm$ and
$\beta\in A^*(T\times Y)$, we have:

\item $\beta/f^! g_*(a)=\bar G^*(\beta)/\bar f^!(a)$
\item $F_!(\beta)/ g_*(a)=\bar F_!\bar G^*(\beta)/a.$
\end{enumerate}
\end{cor}
\Proof
All these relations may be easily obtained using the base-change property.
We illustrate it proving the first relation:
\begin{eqnarray}
f_*(\alpha\cap f^! g_*(a))&=&f_*(\alpha\cap \bar g_*\bar f^!(a))\nonumber\\
&=&f_*\bar g_*(\bar g^*(\alpha)\cap \bar f^!(a)) = g_* \bar f_*
(\bar g^*(\alpha)\cap \bar f^!(a)).
\end{eqnarray}

\Endproof

\begin{property}[\bf Gysin exact sequence]
Let $i\colon Y\mono X$ be a closed embedding and $j\colon X-Y\mono X$ the
corresponding open inclusion. Then, the sequence
\begin{equation*}
A_*(X-Y)\rra {j_*} A_*(X)\rra {i^!} A_*(Y).
\end{equation*}
is exact.
\end{property}

The following  lemma is a ``dualization" of  ``Useful Lemma 1.4.2" from~\cite{Pa}.
\begin{lemma}[\bf Homological useful lemma]\label{useful}
Consider the following diagram with transversal square
$$ \xymatrix{
&X-Y\ar[d]_{k^1}\\
V\ar[r]_{k^0}\ar@/^0.7pc/[d]^p
&X\\
W\ar[r]_{j}\ar[u]^{i} &Y,\ar@{^(->}[u]^{q}}
$$
where $p$ is projective, $q$ is a closed embedding, $X-Y$ is the
open complement of $Y$ in $X$, $k_1$ is the corresponding open
embedding, $pi=\id$, and the map $j$ induces an isomorphism in
homology. Then $\Im k^0_*+\Im k^1_*=A_*(X)$.
\end{lemma}
\Proof Let $x\in A_*(X)$.
Since the map $j_*$ is an isomorphism and $i^!p^!=\id$, we can,
using the transversal base-change property,  lift $x$ up to $\bar
x=p^!(j_*)^{-1}q^!(x)\in A_*(V)$, such that $q^!k_*^0(\bar
x)=q^!(x)$. Then, the Gysin exact sequence implies that
$k_*^0(\bar x)-x\in \Im k^1_*$.
\Endproof

\begin{property}[\bf Projective Bundle Theorem (PBT)]
\label{PBT2}
First, we\newline  should introduce the notion
of an Euler class. For a line bundle $\mathcal L$ over $X$ we set
$e(\mathcal L)\overset{\text{def}}=z^*z_!(1)$,  where $z: X \to
\mathcal L$ is the zero-section (see~\cite[1.1.4]{Pa} for details).

For $X\in\Sm$ and a rank $n$ vector bundle $\E\overset p \to X$
set $\zeta=e(\mathcal O_\E(1))\in A^*(\P(\E))$. Then the
map
\begin{equation*}
\oplus_{i=0}^{n-1}\psi_i \colon A_{*}(\P(\E))\overset\simeq\to
\bigoplus_{i=0}^{r-1}A_{*}(X),
\end{equation*}
where $\psi_i=p_*\circ(\zeta^{\cup i}\cap -)$, is an isomorphism.
\end{property}

\end{section}

\begin{section}{Chern element and homothety involution}
\label{aux2}
 Let $T=\P^1_k$ and take the point $0$ as a
distinguished point. Let $\mathcal A$ be a commutative ring
symmetric $T$-spectrum. For $\lambda \in k^{\times}$
consider a map $\lambda \colon T \to T$ sending $[x\colon y]$ to $[\lambda x\colon y]$
(preserving the distinguished point). It defines an involution
$\lambda^*$ on the sphere $T$-spectrum $\textrm{S}$~\cite{Ya}. For any
space $X$ it defines an involution on the cohomology theory
$A^{*,*}$ as follows:
$$
\epsilon(\lambda)^*=\Sigma_T^{-1}\lambda^*\Sigma_T\colon A^{*,*}(X) \to A^{*,*}(X),
$$
where $\Sigma_T\colon A^{*,*}(X) \to A^{*+2,*+1}(X)$ is the
$T$-suspension isomorphism and $\Sigma_T^{-1}$ is its inverse. Set
\begin{eqnarray}
\label{TheInvolution}
\epsilon=\epsilon(-1)^*.
\end{eqnarray}
The following lemma shows that $\epsilon=\id$ for a certain class of $T$-spectra.

\begin{lemma}
\label{ChernAndEpsilon}
If $\mathcal A$ is equipped with a Chern element
$\gamma \in A^{2,1}(\P^{\infty})$
then
$\epsilon=\id$.
\end{lemma}

\Proof We show that for any $\lambda \in k^{\times}$ one has
$\Sigma_T^{-1}\lambda^*\Sigma_T=\id$. Let $i: \P^1 \to \P^2$ be a
linear embedding. For $X \in \Sm$ consider the projection $p: \P^1
\times X\to \P^1$. The map
$$
(i\times \id_X)^*\colon A^{*,*}(\P^2 \times X) \to A^{*,*}(\P^1 \times X)
$$
is surjective since $1$ and $p^*(\gamma\mid_{\P^1})$ is a free base of the
bigraded $A^{*,*}(X)$-module $A^{*,*}(\P^1 \times X)$.
Now the proof completes as in
\cite[Lemma 1.6]{HY}

\Endproof

\end{section}

\bigskip
\begin{tabbing}
\noindent\textsc{Ivan Panin}\hskip 2cm \={\tt panin@pdmi.ras.ru}\\
\noindent\textsc{Serge Yagunov}\> {\tt yagunov@pdmi.ras.ru}
\end{tabbing}

\noindent\textsc{ Steklov Mathematical Institute (St.Petersburg)}

\noindent\textsc{ Fontanka 27, 191023 St.Petersburg, Russia.}

\end{document}